\documentclass[fleqn]{mat01}
\usepackage{times,mathtimy,amssymb}

\newcommand{\R}{\mbox{$\mathbb{R}$}}
\newcommand{\N}{\mbox{$\mathbb{N}$}}

\def\lin{\mbox{\rm lin}}
\newtheorem{coro}[defin]{\rm COROLLARY}
\begin{document}
\setcounter{page}{253}
\firstpage{253}

\title{Superstability of the generalized orthogonality equation\break on restricted domains}

\author{SOON-MO JUNG and PRASANNA K~SAHOO$^{*}$}

\address{Mathematics Section,
College of Science and Technology,
Hong-Ik University,
339-701 Chochiwon, Korea\\
\noindent $^{*}$Department of Mathematics,
University of Louisville,
Louisville, KY 40292, USA\\
\noindent E-mail: smjung@wow.hongik.ac.kr; sahoo@louisville.edu}

\markboth{Soon-Mo Jung and Prasanna K~Sahoo}{Generalized orthogonality equation on restricted domains}

\volume{114}

\mon{August}

\parts{3}

\Date{MS received 7 August 2003; revised 26 December 2003}

\begin{abstract}
Chmieli\'{n}ski has proved in the paper [4]
the superstability of the generalized orthogonality equation
$|\langle f(x), f(y) \rangle| = |\langle x, y \rangle|$.
In this paper, we will extend the result of Chmieli\'{n}ski by proving
a theorem:
Let $D_{n}$ be a suitable subset of $\,\R^n$.
If a function $f\hbox{:} \ D_{n} \to \R^n$ satisfies the inequality
$||\langle f(x), f(y) \rangle| - |\langle x, y
    \rangle|| ~\leq~ \varphi(x,y)$
for an appropriate control function $\varphi(x,y)$ and for all
$x, y \in D_{n}$, then $f$ satisfies the generalized orthogonality
equation for any $x, y \in D_{n}$.
\end{abstract}


\maketitle


\section{Introduction}

In 1931, Wigner introduced in his book [13] the generalized
orthogonality equation
\begin{equation}
|\langle f(x), f(y) \rangle| ~=~ |\langle x, y \rangle|
\label{eq:zero}
\end{equation}
for all $x, y \in E$, where $E$ is an inner product space and
$\langle\cdot\,,\cdot\rangle$ denotes the inner product on $E$. This
functional equation was solved in \cite{1,2,7,9,10} by many mathematicians.

Recently, Chmieli\'{n}ski [4] proved that the generalized
orthogonality equation is superstable when the relevant functions
belong to the class of functions $f\hbox{:}\ \R^{n} \to \R^n$.
If a function $f\hbox{:}\ \R^n \to \R^n$ ($n \geq 2$) satisfies
the functional inequality
\begin{equation*}
||\langle f(x), f(y) \rangle| - |\langle x, y \rangle|| \leq \varepsilon
\end{equation*}
for some $\varepsilon \geq 0$ and for all $x, y \in \R^n$,
then $f$ is a solution of the generalized orthogonality equation
(\ref{eq:zero}).

We will refer the reader to \cite{3,6,8,12} for detailed definitions of
stability and superstability of functional equations.

By using ideas of Skof and Rassias
\cite{8,11}, and by following the methods of Chmieli\'{n}ski [4,5]
mainly, we will extend the result of Chmieli\'{n}ski by considering
the case when the domain of $f$ is restricted and
by substituting an appropriate control function $\varphi(x,y)$ for
$\varepsilon$ in the relevant inequality as well.

Throughout this paper, let \hbox{$c > 0 \ (c \neq 1)$} and $d > 0$ be
constants and let $n \geq 2$ be a fixed natural number.
By $\N$, $\N_{0}$ and $\R$ we denote the set of
positive integers, of non-negative integers and of real numbers,
respectively.
We will also use the notation $\lin\{ x_{1},\ldots,x_{k} \}$ to denote
the subspace of $\R^n$ spanned by $x_{1},\ldots,x_{k} \in \R^n$.
Let us define a subset $D_{n}$ of $\R^n$ by
\begin{equation*}
D_{n}:= \begin{cases}
                  \{ x \in \R^n: \| x \| \geq d \}, &\mbox{for}\ \  0 < c < 1,\\[.2pc]
                  \{ x \in \R^n: \| x \| < d \},    &\mbox{for}\ \  c > 1,
    \end{cases}
\end{equation*}
where we denote by $\|\!\cdot\!\|$ the usual norm on $\R^n$
defined by
\begin{equation*}
\| x \| ~:=~ \sqrt{\langle x, x \rangle\,}
\end{equation*}
with the usual inner product $\langle\cdot\,,\cdot\rangle$ defined by
\begin{equation*}
\langle x, y \rangle ~:=~ x_{1}y_{1} + x_{2}y_{2} + \cdots +
    x_{n}y_{n}
\end{equation*}
for all points $x = (x_{1},\ldots,x_{n})$ and $y = (y_{1},\ldots,y_{n})$ of $\R^n$.

Suppose $\varphi\hbox{:}\ \R^n \times \R^n \to [0,\infty)$ is a
symmetric function which satisfies the following conditions:
\begin{itemize}
\item[(i)]   There exists a function $\phi\hbox{:}\ [0,\infty) \times [0,\infty) \to [0,\infty)$
                such that $\varphi(x,y) = \phi( \|x\|, \|y\| )$  for all $x, y \in \R^n$.
\item[(ii)]  For all $x, y \in \R^n$,
        \begin{equation*}
        \hskip -1.7pc \frac{1}{|\lambda|} \, \varphi(\lambda x, y) =
                   O\!\left( -\frac{\ln c}{\ln |\lambda|} \right)
        \end{equation*}
                 either as $| \lambda | \to \infty$ (for $0 < c < 1$) or as $| \lambda | \to 0$ (for $c>1$).
\item[(iii)] If both $|\lambda|$ and $|\mu|$ are different from 1, then for all $x, y \in \R^n$,
                 \begin{equation*}
         \hskip -1.7pc \frac{1}{|\lambda\mu|} \, \varphi(\lambda x, \mu y)
                   = O\!\left( \left| \frac{\ln c}{\ln\!|\lambda|} \,
            \frac{\ln c}{\ln\!|\mu|} \right|\,\right)
         \end{equation*}
                 either as $| \lambda\mu | \to \infty$ (for $0 < c < 1$) or as $|\lambda\mu| \to 0$ (for $c>1$).
\end{itemize}


\section{Preliminaries}

We begin by introducing a lemma of \cite{4} which
turns out to be very useful to prove Lemma 4 below.


\renewcommand\thedefin{\arabic{defin}}

\begin{theorem}[\!]
Let $\varepsilon \geq 0$ be given.
For each $\eta > 0$ there exists $k_{0} \in \N$ such that if $a${\rm ,}
$u_{1}${\rm ,} $u_{2},\ldots,u_{n-1} \in \R^{n}\!\setminus\!\{ 0 \}$
satisfy the conditions
\begin{align*}
&1 - \frac{\varepsilon}{k^{2}} \leq \| u_{i} \|^2 ~\leq~
    1 + \frac{\varepsilon}{k^{2}} \quad (i=1,2,\ldots,n-1),\\[.2pc]
&|\langle u_{i}, u_{j} \rangle| \leq
    \frac{\varepsilon}{k^{2}}  \quad (i,j=1,2,\ldots,n-1;\,i \neq j),\\[.2pc]
&|\langle a, u_{i} \rangle| \leq \frac{\varepsilon}{k}
    \quad(i=1,2,\ldots,n-1){\rm ,}
\end{align*}
for any $k \geq k_{0}${\rm ;} then
\begin{itemize}
\leftskip 1.3pc
\item[{\rm (a)}] $u_{1},\ldots,u_{n-1}$ are linearly independent{\rm ;}
\item[{\rm (b)}] $|\cos A(a,\ell)| \geq 1 - \eta$, where $\ell$
               denotes the line in $\R^n$ which is the
               orthogonal complement of\ $\lin\{ u_{1},\ldots,u_{n-1} \}$
               and $A(\cdot\,,\cdot)$ stands for the angle.
\end{itemize}
\end{theorem}

In the following five lemmas, we will modify the statements of
Proposition 1 in \cite{4} and later apply them to the
proof of our main result.\pagebreak

In the following lemmas and theorems of this section, we assume
that the function $f\hbox{:}\ D_{n} \to \R^n$ satisfies the inequality
\begin{equation}
||\langle f(x), f(y) \rangle| - |\langle x, y \rangle|| \leq \varphi(x,y)
\label{eq:one}
\end{equation}
for all $x, y \in D_{n}$ if there is no specification for $f$.

It is enough to put $y=x$ in the inequality (\ref{eq:one}) to prove
the following lemma.


\begin{lemma}
The following inequality
\[\| x \|^2 - \varphi(x,x) ~\leq~ \| f(x) \|^2 ~\leq~ \| x \|^2 +
    \varphi(x,x)\]
holds for any $x \in D_{n}$.
\end{lemma}


\begin{lemma}
If $f(x)=0${\rm ,} then $x=0$.
\end{lemma}

\begin{proof}
If $f(x)=0$, then (\ref{eq:one}) implies for each
$y \in D_{n}$ that $|\langle x, y \rangle| \leq \varphi(x,y)$.
By putting $y=\lambda x \in D_{n}$ in the last inequality, we obtain
\begin{equation}
\| x \|^2 ~\leq~ \frac{1}{| \lambda |}\,\varphi(x,\lambda x).
\label{eq:two}
\end{equation}

If $0<c<1$ and if we take the limit in (\ref{eq:two}) as
$| \lambda | \to \infty$, then (ii) implies $x=0$ which is
impossible because $\| x \| \geq d > 0$.
Thus, if $0<c<1$, then $f(x)\neq 0$ for every $x \in D_{n}$.
When $c>1$, we can take the limit in (\ref{eq:two}) as
$| \lambda | \to 0$.
Then, (ii) and (\ref{eq:two}) yield $x=0$.\hfill$\Box$
\end{proof}


\begin{lemma}
For each $x \in D_{n}\!\setminus\!\{ 0 \}$
there exists a function $\mu_{x}\hbox{:}\ \R \to \R$ such that
$f(\lambda x) = \mu_{x}(\lambda) f(x)$ for all
$\lambda\in\R\!\setminus\!\{ 0 \}$ with $\lambda x \in D_{n}$
and also such that
\begin{equation*}
\frac{| \mu_{x}(\lambda) |}{| \lambda |} \not\to 0 \;\;
    \begin{cases}
    \mbox{as}\; |\lambda| \to \infty &(\mbox{for}\ \ \  0<c<1)\\[2pt]
        \mbox{as}\; |\lambda| \to 0      &(\mbox{for}\ \ \ c>1)
    \end{cases}.
\end{equation*}
\end{lemma}

\begin{proof}
Assume that $x\in D_{n}\!\setminus\!\{ 0 \}$ and
$\lambda\neq 0$ are given with $\lambda x \in D_{n}$.
If $f(x)$ and $f(\lambda x)$ were linearly independent,
then we could select some $\omega>0$ such that
\begin{equation}
| \cos A(f(x), f(\lambda x)) | ~=~ 1-\omega,  \label{eq:three}
\end{equation}
where $A(\cdot\,,\cdot)$ stands for the angle.

Since $x \neq 0$ is assumed, we can choose an orthogonal basis
$\{ x, v_{1},\ldots,v_{n-1} \}$ for $\R^n$ with
$\| v_1 \| = \cdots = \| v_{n-1} \| = 1$.
For any $i,j\in\{ 1,\ldots,n-1 \}$ with $i \neq j$ and for any
$k \in \N$ we have
\begin{equation}
\langle c^{-k}v_{i}, x \rangle = \langle c^{-k}v_{i}, \lambda x
\rangle = \langle c^{-k}v_{i}, c^{-k}v_j \rangle = 0.
\label{eq:four}
\end{equation}

By using simple notations given by
\begin{equation*}
a := f(x), \;\;\; a' := f(\lambda x) \quad\mbox{ and }\quad
    u_i := c^k f(c^{-k}v_i ),
\end{equation*}
we get from (\ref{eq:one}), (\ref{eq:four}), (ii) and (iii) that
\begin{align*}
|\langle u_{i}, a \rangle| &\leq c^k \varphi(c^{-k}v_{i},x)=O\!\left( \frac{1}{k} \right),\\[2pt]
|\langle u_{i}, a' \rangle| &\leq c^k \varphi(c^{-k}v_{i},\lambda x)=O\!\left( \frac{1}{k} \right),\\[2pt]
|\langle u_i, u_j \rangle| &\leq c^{2k} \varphi(c^{-k}v_{i},c^{-k}v_{j})=O\!\left( \frac{1}{k^2} \right)
\end{align*}
for all $i,j \in \{ 1,\ldots,n-1 \}$ with $i\neq j$ and for any
sufficiently large $k\in\N$ (such that $c^{-k}v_i \in D_{n}$
for $i=1,\ldots,n-1$).
Moreover, it follows from Lemma 2 that
\begin{equation*}
1 - c^{2k} \varphi(c^{-k}v_{i}, c^{-k}v_{i}) ~\leq~ \| u_i \|^2 \leq 1 + c^{2k} \varphi(c^{-k}v_{i}, c^{-k}v_{i}).
\end{equation*}

At this point, we apply Theorem 1.
First, denote by $\ell$ the one-dimensional orthogonal complement
of the subspace $\lin\{ u_{1},\ldots,u_{n-1} \}$.
According to Theorem 1, (ii) and (iii), we can choose a
sufficiently large integer $k$ in order that $| \cos A(a,\ell) |$
and $| \cos A(a',\ell) |$ are arbitrarily close to 1. This fact
means that $| \cos A(a,a') | = | \cos A(f(x),f(\lambda x)) |$ is
really 1, which is contrary to our assumption (\ref{eq:three}).
Therefore, $f(x)$ and $f(\lambda x)$ have to be linearly
dependent.

According to Lemma 3, $f(x)\neq 0$ and $f(\lambda x)\neq 0$
because $x\in D_{n}\!\setminus\!\{ 0 \}$ and $\lambda\neq 0$ with
$\lambda x \in D_{n}$.
Thus, we can choose $\mu_{x}(\lambda) \in \R$ such that
$f(\lambda x) = \mu_{x}(\lambda) f(x)$ for all
$\lambda\neq 0$ with $\lambda x \in D_{n}$.

Hence, with $y = \lambda x$, (\ref{eq:one}) yields
\begin{equation}
\| x \|^2 - \frac{1}{| \lambda |}\,\varphi(x,\lambda x) ~\leq~
\frac{| \mu_{x}(\lambda) |}{| \lambda |}\,\| f(x) \|^2    ~\leq~
\| x \|^2 + \frac{1}{| \lambda |}\,\varphi(x,\lambda x)
\label{eq:five}
\end{equation}
for any $x\in D_{n}\!\setminus\!\{ 0 \}$ and $\lambda\neq 0$ with
$\lambda x \in D_{n}$.

When $0<c<1$, by taking the limit in (\ref{eq:five}) as
$| \lambda | \to \infty$, (\ref{eq:five}) and (ii) imply that
\begin{equation*}
\frac{| \mu_{x}(\lambda) |}{| \lambda |} \not\to 0
    \quad\mbox{as} \ |\lambda| \to \infty.
\end{equation*}
Similarly, when $c>1$, we can take the limit in (\ref{eq:five}) by
letting $| \lambda | \to 0$ and use (\ref{eq:five}) and (ii) to
obtain
\begin{equation*}
\frac{| \mu_{x}(\lambda) |}{| \lambda |} \not\to 0\quad\mbox{ as } | \lambda | \to 0,
\end{equation*}
which completes the proof.\hfill$\Box$
\end{proof}


\begin{lemma}
For any $x \in D_{n}\!\setminus\!\{ 0 \}$ and
$y \in D_{n}${\rm ,} it holds that $\langle x, y \rangle = 0$ if and only if
$\langle f(x), f(y) \rangle = 0$.
\end{lemma}

\begin{proof}
Let $x \in D_{n}\!\setminus\!\{ 0 \}$, $y \in D_{n}$
and $\lambda \neq 0$ be given with $\lambda x \in D_{n}$.
If $\langle x, y \rangle = 0$,
then $\langle \lambda x, y \rangle = 0$.
In this case, it follows from (\ref{eq:one}) that
\begin{equation}
|\langle f(\lambda x), f(y) \rangle| \leq \varphi(\lambda x, y).
\label{eq:seven}
\end{equation}

On account of Lemma 4, there exists a function
$\mu_{x}\hbox{:}\ \R \to \R$ such that
\begin{equation}
f(\lambda x) ~=~ \mu_{x}(\lambda) f(x)  \label{eq:eight}
\end{equation}
for all $x \in D_{n}\!\setminus\!\{ 0 \}$ and $\lambda \neq 0$ with
$\lambda x \in D_{n}$.
Moreover, there is a constant $\alpha > 0$ and a strictly
increasing (or decreasing) positive sequence $(\lambda_{k})$ with
\begin{equation*}
\begin{cases}
       \lambda_{k} \to \infty,  &\mbox{ for } \ \ \ 0<c<1,\\[2pt]
       \lambda_{k} \to 0,       &\mbox{ for } \ \ \ c>1,
\end{cases}
\end{equation*}
such that
\begin{equation}
\frac{| \mu_{x}(\lambda_{k}) |}{\lambda_{k}} ~\geq~ \alpha
\label{eq:nine}
\end{equation}
for every $k \in \N$.
By (\ref{eq:seven})--(\ref{eq:nine}) and (ii),
we have
\begin{align}
|\langle f(x), f(y) \rangle|
     &\leq  \frac{1}{| \mu_{x}(\lambda_{k}) |}\,
              \varphi(\lambda_{k}x,y)\nonumber\\[2pt]
     &\leq  \frac{1}{\alpha\lambda_{k}}\,
               \varphi(\lambda_{k}x,y)\nonumber\\[2pt]
     &\to   0, \quad \mbox{ as }\; \ k \to \infty. \label{eq:ten}
\end{align}

Suppose $x \in D_{n}\!\setminus\!\{ 0 \}$ and $y \in D_{n}$ are
given with $\langle f(x), f(y) \rangle = 0$.
For each $\lambda > 0$ with $\lambda x \in D_{n}$, Lemma 4 gives
\begin{equation*}
|\langle f(\lambda x), f(y) \rangle| =
    | \mu_{x}(\lambda) | \, |\langle f(x), f(y) \rangle| = 0.
\end{equation*}
Hence, it follows from (\ref{eq:one}) and (ii) that
\begin{align*}
|\langle x, y \rangle|
     &\leq  \frac{1}{\lambda}\,\varphi(\lambda x, y)\\[2pt]
     &\to   0\;\begin{cases}
                 \mbox{as}\ \  \lambda \to \infty, &(\mbox{for}\ \  0<c<1),\\[2pt]
                 \mbox{as}\ \  \lambda \to 0,      &(\mbox{for}\ \  c>1),
              \end{cases}
\end{align*}
which finishes our proof.\hfill$\Box$
\end{proof}

In the following lemma, we will prove the converse of Lemma 3, i.e.,
$f(0) = 0$ under an essential condition that the range of $f$ is
a finite-dimensional space.


\begin{lemma}
It holds that $f(0) = 0$.
\end{lemma}

\begin{proof}
We need to consider the case $c>1$ only because
$D_{n}$ does not contain 0 for the other case $0<c<1$.
For any $x \in D_{n}\!\setminus\!\{ 0 \}$, by putting $y=0$,
Lemma 5 gives
\begin{equation}
\langle f(x), f(0)\rangle=0.  \label{eq:eleven}
\end{equation}

Let $\{ x_{1}, \ldots, x_{n} \}$ be an orthogonal basis for
$\R^n$ with $x_{i} \in D_{n}\!\setminus\!\{ 0 \}$.
Then, Lemma 3 implies $f(x_i) \neq 0$ for $i=1,\ldots,n$.
Furthermore, Lemma 5 implies that
$\langle f(x_i), f(x_j) \rangle = 0$ for any
$i,j \in \{ 1,\ldots,n \}$ with $i \neq j$, i.e.,
$\{ f(x_1),\ldots,f(x_n) \}$ is another orthogonal basis for
$\R^n$.

Therefore, it follows from (\ref{eq:eleven}) that
\begin{equation*}
\langle f(x_i), f(0) \rangle=0
\end{equation*}
for $i=1,\ldots,n$, and this relation implies $f(0)=0$.\hfill$\Box$
\end{proof}

By using ideas from Propositions 1 and 2 of [5], we can
prove the following theorem.


\begin{theorem}[\!]
Assume that a function $f\!:D_n \to \R^n$ satisfies the
functional inequality
\begin{equation*}
||\langle f(x), f(y) \rangle| - |\langle x, y \rangle||\leq \varphi(x,y)
\end{equation*}
for all $x, y \in D_n$.
We have for $2 \leq k \leq n${\rm :}
\begin{itemize}
\leftskip 1.3pc
\item[{\rm (a)}] $x_{1},\ldots,x_{k} \in D_{n}$ are linearly independent
               if and only if $f(x_{1}),\ldots,f(x_{k})$ are linearly
               independent{\rm ;}
\item[{\rm (b)}] Let ${\cal P}$ be a $k$-dimensional subspace of
               $\,\R^n$.
               Then $f$ transforms ${\cal P} \cap D_{n}$ into
               $k$-dimensional subspace ${\cal P'}$ of $\,\R^n$
               spanned by the images of the elements of an arbitrary
               basis ${\cal B}$ of ${\cal P}$ with
               ${\cal B} \subset {\cal P} \cap D_{n}$.
\end{itemize}
\end{theorem}

\begin{proof}$\left.\right.$\vspace{.5pc}

\noindent (a) Let $x_{1},\ldots,x_{k} \in D_{n}$ be linearly
independent and suppose that $f(x_{1})$, \ldots, $f(x_{k})$ are
linearly dependent.
Then, we can choose $\lambda_{2},\ldots,\lambda_{k} \in \R$ such
that
\begin{equation}
f(x_{1}) = \lambda_{2} f(x_{2}) + \cdots + \lambda_{k} f(x_{k}).
\label{eq:twelve}
\end{equation}

Let $x \in \lin \{ x_{1},\ldots,x_{k} \} \cap D_{n}$ be chosen with
$x \neq 0$ and $\langle x, x_{i} \rangle = 0$ for
$i=2,\ldots,k$.
According to Lemma 5, it holds that $\langle f(x), f(x_{i}) \rangle = 0$
for $i=2,\ldots,k$, and hence (\ref{eq:twelve}) implies
$\langle f(x), f(x_{1}) \rangle = 0$.
By Lemma 5 again, $\langle x, x_{i} \rangle = 0$ for $i=1,\ldots,k$,
and hence $x \not\in \lin \{ x_{1},\ldots,x_{k} \}$,
a contradiction.

In this paper, the converse of the above statement will not be used.
But here we will introduce its proof for completion.
Let $f(x_{1}),\ldots,f(x_{k})$ be linearly independent and
$x_{1},\ldots,x_{k}$ be linearly dependent.
Then, there are real numbers $\lambda_{2},\ldots,\lambda_{k}$ such
that
\begin{equation}
x_{1} = \lambda_{2}x_{2} + \cdots + \lambda_{k}x_{k}.
\label{eq:add1}
\end{equation}

Choose $y \in \lin \{ f(x_{1}),\ldots,f(x_{k}) \} \cap f(D_{n})$ with
$y \neq 0$ and $\langle y, f(x_{i}) \rangle = 0$ for
$i=2,\ldots,k$.
There exists an $x \in D_{n}\!\setminus\!\{ 0 \}$ with $y = f(x)$.
Due to Lemma 5, we have $\langle x, x_{i} \rangle = 0$
for $i=2,\ldots,k$, and (\ref{eq:add1}) means
$\langle x, x_{1} \rangle = 0$.
Using Lemma 5 again, we obtain
$\langle f(x), f(x_{i}) \rangle = \langle y, f(x_{i}) \rangle = 0$
for $i=1,\ldots,k$.
This implies that $y \not\in \lin \{ f(x_{1}),\ldots,f(x_{k}) \}$
which leads to a contradiction.

\noindent (b) Let $\{ x_{1},\ldots,x_{k} \} \subset D_{n}$ be an orthogonal basis
for a $k$-dimensional subspace ${\cal P}$ of $\R^n$ and let
$\{ x_{1},\ldots,x_{k},x_{k+1},\ldots,x_{n} \}$ $\subset$ $D_{n}$ be an
orthogonal basis for $\R^n$.
On account of Lemmas 3 and 5,
$\{$ $\!\!f(x_{1})$, \ldots, $f(x_{n})$ $\!\!\}$
is also an orthogonal basis for $\R^n$.
Thus, for any $x \in {\cal P} \cap D_{n}$ there exist
$\lambda_{1}$, $\ldots$, $\lambda_{k}$, $\xi_{1}$, \ldots, $\xi_{n}$ $\in$
$\R$ such that
\begin{equation}
x = \lambda_{1} x_{1} + \cdots + \lambda_{k} x_{k}
\quad\mbox{ and }\quad
f(x) = \xi_{1} f(x_{1}) + \cdots + \xi_{n} f(x_{n}).
\label{eq:thirteen}
\end{equation}

Since $\langle x, x_{i} \rangle = 0$ for $i=k+1,\ldots,n$, Lemma 5
implies that $\langle f(x), f(x_{i}) \rangle = 0$ for $i=k+1,\ldots,n$.
Hence, it follows from (\ref{eq:thirteen}) that
\begin{equation*}
\langle f(x), f(x_{i}) \rangle ~=~ \xi_{i}\,\|f(x_{i})\|^2=0
\end{equation*}
for $i=k+1,\ldots,n$, and we have $\xi_{k+1} = \cdots = \xi_{n} = 0$.
Therefore, we conclude that
\begin{equation*}
f(x) = \xi_{1} f(x_{1}) + \cdots + \xi_{k} f(x_{k})
\end{equation*}
or
\begin{equation*}
f({\cal P} \cap D_{n}) \subset \lin \{ f(x_{1}),\ldots,f(x_{k}) \}.
\end{equation*}

If $\{ y_{1},\ldots,y_{k} \} \subset D_{n}$ is a basis for ${\cal P}$,
it then follows from (a) that $\{ f(y_{1}),\ldots,f(y_{k}) \}$ is a
basis for $\lin \{ f(x_{1}),\ldots,f(x_{k}) \}$, and this completes
the proof.\hfill$\Box$
\end{proof}

In the following lemma, we will modify Lemma 3 of
[4] in order to be applicable to our case.


\begin{lemma}
It holds that
\begin{equation*}
\lim_{k \to \infty} c^{2k} \, \| f(c^{-k}x) \| \, \| f(c^{-k}y) \|
    = \| x \| \, \| y \|
\end{equation*}
for all $x,y \in D_{n}$.
\end{lemma}

\begin{proof}
By Lemma 2, we get
\begin{equation*}
\begin{array}{l}
       \sqrt{\|x\|^2 - c^{2k} \varphi(c^{-k}x, c^{-k}x)\,} \,
       \sqrt{\|y\|^2 - c^{2k} \varphi(c^{-k}y, c^{-k}y)\,}
       \vspace{3mm}\\
       \leq~ c^{2k} \, \| f(c^{-k}x) \| \, \| f(c^{-k}y) \|
      \vspace{3mm}\\
      \leq \sqrt{\|x\|^2 + c^{2k} \varphi(c^{-k}x, c^{-k}x)\,} \,
       \sqrt{\|y\|^2 + c^{2k} \varphi(c^{-k}y, c^{-k}y)\,},
    \end{array}
\end{equation*}
and (iii) gives the validity of our assertion.
\hfill$\Box$
\end{proof}

\section{Main results}

We know that $D_{n}$ is a subset of $\R^n$ defined by
\begin{equation*}
D_{n} := \begin{cases}
                  \{ x \in \R^n\hbox{:}\ \| x \| \geq d \}, &\mbox{for}\ \  0 < c < 1,\\[2pt]
                  \{ x \in \R^n\hbox{:}\ \| x \| < d \},    &\mbox{for}\ \  c > 1,
    \end{cases}
\end{equation*}
for given positive numbers $c \neq 1$ and $d > 0$.
The function $\varphi\hbox{:}\ \R^n \times \R^n \to [0,\infty)$ was
defined as a symmetric function which satisfies the following
conditions:
\begin{itemize}
\item[(i)]   There exists a function
                $\phi\hbox{:}\ [0,\infty) \times [0,\infty) \to [0,\infty)$
                such that $\varphi(x,y) = \phi( \|x\|, \|y\| )$
                for all $x, y \in \R^n$.
\item[(ii)]  For all $x, y \in \R^n$,
\begin{equation*}
\hskip -1.7pc \frac{1}{|\lambda|} \, \varphi(\lambda x, y) = O\!\left( -\frac{\ln c}{\ln\!|\lambda|} \right)
\end{equation*}
                 either as $| \lambda | \to \infty$ (for $0 < c < 1$)
                 or as $| \lambda | \to 0$ (for $c>1$).
\item[(iii)] If both of $|\lambda|$ and $|\mu|$ are different from
                 1, then for all $x, y \in \R^n$,
\begin{equation*}
\hskip -1.7pc \frac{1}{|\lambda\mu|} \, \varphi(\lambda x, \mu y)=O\!\left( \left| \frac{\ln c}{\ln\!|\lambda|} \,
                                 \frac{\ln c}{\ln\!|\mu|} \right|\,\right)
\end{equation*}
                 either as $| \lambda\mu | \to \infty$ (for $0 < c < 1$)
                 or as $|\lambda\mu| \to 0$ (for $c>1$).
\end{itemize}

As assumed in the previous section,
throughout this section also, let the function
$f\hbox{:}\ D_{n} \to \R^n$ satisfy the functional inequality
(\ref{eq:one}) for all $x, y \in D_{n}$ if there is no
specification for $f$.


\begin{lemma}
It holds that
\begin{equation*}
|\cos A(f(x), f(y))|=|\cos A(x,y)|
\end{equation*}
for any $x$ and $y$ in $D_{n}\!\setminus\!\{ 0 \}$.
\end{lemma}

\begin{proof}
By making use of Lemmas 3 and 4, it is easy to see
\begin{equation}
|\cos A(f(x), f(y))|=|\cos A(f(c^{-k}x), f(c^{-k}y))|
\label{eq:forteen}
\end{equation}
for all $x, y \in D_{n}\!\setminus\!\{ 0 \}$ and any $k \in \N$.

If we replace $x, y$ in (\ref{eq:one}) by $c^{-k}x$ and $c^{-k}y$,
respectively, and if we divide the resulting inequalities by
$c^{-2k}$, then
\begin{align*}
&\| x \| \, \| y \| \, | \cos A(x,y) |
- c^{2k}\varphi(c^{-k}x, c^{-k}y)\\[2pt]
&\quad\leq c^{2k} \, \| f(c^{-k}x) \| \, \| f(c^{-k}y) \| \,
| \cos A(f(c^{-k}x),  f(c^{-k}y)) |\\[2pt]
&\quad\leq \| x \| \, \| y \| \, | \cos A(x,y) |
+ c^{2k}\varphi(c^{-k}x, c^{-k}y).
\end{align*}

Taking the limit as $k \to \infty$ in the above inequalities and
using (iii), (\ref{eq:forteen}) and Lemma 8, we obtain
\begin{equation*}
\| x \| \, \| y \| \, | \cos A(x,y) | =
    \| x \| \, \| y \| \, | \cos A(f(x), f(y)) |,
\end{equation*}
which ends the proof.\hfill$\Box$
\end{proof}

We now define an integer $k_{0} \in \N_{0}$\  by
\begin{equation*}
k_{0} := \min \{ k \in \N_{0}\hbox{:}\; c^{-k}e_{i} \in D_{n}
                       \quad\mbox{for all }\ \  i=1,\ldots,n \},
\end{equation*}
and let
\begin{equation*}
e'_{i} ~:=~ c^{-k_{0}}e_{i}
\end{equation*}
for $i=1,\ldots,n$, where $\{ e_{1},\ldots,e_{n} \}$ is the canonical basis for $\R^{n}$.


\begin{lemma}
There exists an orthogonal automorphism
$\psi\hbox{:}\ \R^n \to \R^n$ such that
\begin{itemize}
\leftskip 1.3pc
\item[{\rm (a)}] the composition $f':= \psi \circ f$ satisfies the
               inequality $(\ref{eq:one})$ for all $x, y \in D_{n}${\rm ;}
\item[{\rm (b)}] every element of $\,\{ e'_{1},\ldots,e'_{n} \}$ is an
               eigenvector of $f'$, i.e.{\rm ,}
               \begin{equation*}
        f'(e'_{i}) ~=~ \lambda_{i} e'_{i},
           \end{equation*}
               where $\lambda_{i}$ is a constant with
               $0 < \lambda_{i} \leq
                \sqrt{1+c^{2k_{0}}\varphi(e'_{i},e'_{i})\,}\,$ for
               $i = 1,\ldots,n$.
\end{itemize}
\end{lemma}

\begin{proof}$\left.\right.$\vspace{.5pc}

\noindent (a) By Lemmas 3 and  5,
$\{ f(e'_{1}),\ldots,f(e'_{n}) \}$ is an orthogonal basis for
$\R^n$.
We may define an orthogonal automorphism
$\psi\hbox{:}\ \R^n \to \R^n$ by
\begin{equation}
\psi(x) := \lambda_{1} \, \| f(e'_{1}) \| \, e_{1} + \cdots +
               \lambda_{n} \, \| f(e'_{n}) \| \, e_{n}
\label{eq:fifteen}
\end{equation}
for any $x \in \R^n$ expressed by
$x = \lambda_{1} f(e'_{1}) + \cdots + \lambda_{n} f(e'_{n})$.
Since $\psi$ is orthogonal, we have
\begin{equation*}
\langle \psi(f(x)), \psi(f(y)) \rangle =
    \langle f(x), f(y) \rangle
\end{equation*}
for all $x, y \in D_{n}$.
Hence, it is obvious that $f' = \psi \circ f$ satisfies
inequality (\ref{eq:one}) for all $x, y \in D_{n}$.

\noindent (b) By (\ref{eq:fifteen}), we obtain
\begin{equation*}
f'(e'_{i}) = \psi(f(e'_{i})) = \| f(e'_{i}) \| \, e_{i}
    = c^{k_{0}}\,\| f(e'_{i}) \|\, e'_{i}.
\end{equation*}
Further, it follows from Lemmas 2 and 3 that
\begin{equation*}
0 < c^{k_{0}}\,\| f(e'_{i}) \| \leq
    \sqrt{1+c^{2k_{0}}\varphi(e'_{i},e'_{i})\,}
\end{equation*}
for $i=1,\ldots,n$.\hfill$\Box$
\end{proof}

On the basis of Theorem 7(b), we are now ready to deal with
a special case of $n=2$ in the following lemma.


\begin{lemma}
Let a function $f\!:\!D_{2} \to \R^2$
satisfy the inequality $(\ref{eq:one})$ for all $x,y\in D_{2}$.
If $f(e'_{i}) = \lambda_{i} e'_{i}$ for $i=1,2$ with
\begin{equation*}
   0 < \lambda_{i} ~\leq~
    \sqrt{1+c^{2k_{0}}\varphi(e'_{i},e'_{i})\,},
\end{equation*}
then either $f(x) = x${\rm ,} $f(x) = -x${\rm ,} $f(x) = \overline{x}${\rm ,} or
$f(x) = -\overline{x}$ for each $x \in D_{2}${\rm ,} where
$\overline{x} = (x_{1}, -x_{2})$ for $x = (x_{1}, x_{2})$ and
see Lemma {\rm 10} for the $e'_{i}$'s and $k_{0}$.
\end{lemma}

\begin{proof}
According to Lemma 6, it holds $f(0) = 0$.
This means the validity of our assertion for $x = 0$
(if 0 belongs to $D_{2}$).

Now, let $x \in D_{2}\!\setminus\!\{ 0 \}$.
Due to Lemma 9, we have
\begin{equation*}
|\cos A(x,e'_{1})|=|\cos A(f(x),f(e'_{1}))|=|\cos A(f(x),e'_{1})|.
\end{equation*}
This implies that there exists a non-zero real number
$\lambda$ such that either $f(x) = \lambda x$ or
$f(x) = \lambda\overline{x}$.

Let us define
\begin{equation*}
{\cal A} := \{ x \in D_{2}\hbox{:}\ \mbox{there exists a}\
    \lambda\in\R \ \mbox{with}\ f(x) = \lambda x \}
\end{equation*}
and
\begin{equation*}
{\cal B} := \{ x \in D_{2}\hbox{:}\ \mbox{there exists a}\
    \lambda\in\R \ \mbox{with} \ f(x) = \lambda\overline{x} \}.
\end{equation*}

On account of Lemma 4, it is not difficult to see
$D_{2} \cap (\lin\{ e_{1} \} \cup \lin\{ e_{2} \}) \subset
   {\cal A} \cap {\cal B}$.
We set
\begin{equation*}
D_{2}^{*} := D_{2}\!\setminus\!(\lin\{ e_{1} \} \cup\lin\{ e_{2} \}).
\end{equation*}
We assert that either $D_{2}^{*} \subset {\cal A}$ or
$D_{2}^{*} \subset {\cal B}$.
Suppose that there were $x,y \in D_{2}^{*}$ and
$\lambda, \mu \in \R\!\setminus\!\{ 0 \}$ with
\begin{equation*}
f(x)=\lambda x \quad\mbox{and}\quad f(y) = \mu \overline{y}.
\end{equation*}
By Lemma 9, we would have
\begin{equation*}
|\cos A(x,\overline{y})| = |\cos A(f(x),f(y))|=|\cos A(x,y)|.
\end{equation*}
This implies that $x$ or $y$ should belong to
$D_{2} \cap (\lin\{ e_{1} \} \cup \lin\{ e_{2} \})$,
which leads to a contradiction.

From the above fact we can deduce that there exists a function
$\lambda\hbox{:}\ D_{2} \to \R$ such that either
\begin{equation}
f(x)=\lambda(x) x \quad\mbox{ for all }\ x \in D_{2}\label{eq:sixteen}
\end{equation}
or
\begin{equation}
f(x)=\lambda(x) \overline{x} \quad\mbox{ for all}\ x \in D_{2}.
\label{eq:seventeen}
\end{equation}

In (\ref{eq:sixteen}), it follows from (\ref{eq:one})
that
\begin{equation*}
||\langle \lambda(x) x, \lambda(x) x \rangle| -
    |\langle x, x \rangle|| ~\leq~ \varphi(x,x)
\end{equation*}
for any $x \in D_{2}$, and hence
\begin{equation*}
| \lambda(x)^{2} - 1 |\,\|x\|^2 ~\leq~ \varphi(x,x)
\end{equation*}
for $x\in D_{2}$.
If we replace $x$ by $c^{-k}x$ in the last inequality, then we get
\begin{equation*}
| \lambda(c^{-k}x)^2 - 1 |\,\|x\|^2 \leq c^{2k}\varphi(c^{-k}x,c^{-k}x),
\end{equation*}
and if we take the limit as $k \to \infty$, then (iii) means
\begin{equation}
\lim_{k\to\infty} \lambda(c^{-k}x)^2=1
\label{eq:eightteen}
\end{equation}
for any $x \in D_{2}\!\setminus\!\{ 0 \}$.
Choose $x,y \in D_{2}$ with $\langle x, y \rangle \neq 0$ and let
$k \in \N$.
It follows from (\ref{eq:one}) that
\begin{equation*}
||\langle \lambda(x)x, \lambda(c^{-k}y)c^{-k}y \rangle| -
    |\langle x, c^{-k}y \rangle|| ~\leq~ \varphi(x,c^{-k}y).
\end{equation*}
By making use of (ii) and (\ref{eq:eightteen}) and
by taking the limit as $k \to \infty$, we conclude that
$| \lambda(x) | = 1$ for every $x \in D_{2}\!\setminus\!\{ 0 \}$,
i.e.,
\begin{equation*}
f(x)=x \quad\mbox{ or }\quad f(x)=-x,
\end{equation*}
for all $x \in D_{2}$, in view of (\ref{eq:sixteen}) and Lemma 6.

In (\ref{eq:seventeen}), we can analogously obtain
the equality (\ref{eq:eightteen}) for each
$x \in D_{2}\!\setminus\!\{ 0 \}$ because of the fact
$\|\overline{x}\| = \|x\|$.
The fact
$\langle \overline{x}, \overline{y} \rangle = \langle x, y \rangle$
yields $| \lambda(x) | = 1$ for each
$x \in D_{2}\!\setminus\!\{ 0 \}$
and hence (\ref{eq:seventeen}) and Lemma 6 give
\begin{equation*}
f(x)=\overline{x} \quad\mbox{or}\quad f(x) = -\overline{x},
\end{equation*}
for all $x \in D_{2}$, which completes the proof.\hfill$\Box$
\end{proof}

By making use of Lemmas 10 and 11 we can easily prove the
following corollary.
Hence, we omit the proof.


\begin{coro}$\left.\right.$\vspace{.5pc}

\noindent If a function $f\!:\!D_{2} \to \R^2$
satisfies the inequality $(\ref{eq:one})$ for all $x,y \in D_{2}${\rm ,}
then
\begin{equation*}
\| f(x) \| ~=~ \| x \|
\end{equation*}
for every $x$ in $D_{2}$.
\end{coro}

In the following lemma, we will extend the last corollary to the
spaces of higher dimensions.


\begin{lemma}
If a function $f\!:\!D_{n} \to \R^n$
satisfies the inequality $(\ref{eq:one})$ for all $x,y \in D_{n}${\rm ,}
then
\begin{equation*}
\| f(x) \| ~=~ \| x \|
\end{equation*}
for any $x \in D_{n}$.
\end{lemma}

\begin{proof}
Lemma 6 says that $f(0) = 0$, and this means that our
assertion holds true for $x=0$ whenever $0\in D_{n}$ (i.e., in the
case $c > 1$).

We now choose $x, y \in D_{n}\!\setminus\!\{ 0 \}$ with
$\langle x, y \rangle = 0$.
In view of Lemmas 3 and 5, we know that $f(x) \neq 0$,
$f(y) \neq 0$ and $\langle f(x), f(y) \rangle = 0$.
Due to Theorem 7(b) and Lemma 4, we obtain
\begin{equation*}
f\!\left( D_{n} \cap \lin\left\{ \frac{x}{\|x\|},  \frac{y}{\|y\|}
    \right\}\right)\subset \lin\left\{ \frac{f(x)}{\|f(x)\|},
    \frac{f(y)}{\|f(y)\|} \right\}.
\end{equation*}
This means that for each pair $(\lambda_{1}, \lambda_{2})$ of real
numbers satisfying
\begin{equation}
\lambda_{1}\frac{x}{\|x\|} + \lambda_{2}\frac{y}{\|y\|} ~\in~
D_{n},  \label{eq:nineteen}
\end{equation}
there exists a unique pair $(\mu_{1}, \mu_{2})$ of real numbers
such that
\begin{equation}
f\!\left( \lambda_{1}\frac{x}{\|x\|} + \lambda_{2}\frac{y}{\|y\|}
\right) ~=~ \mu_{1}\frac{f(x)}{\|f(x)\|} +
\mu_{2}\frac{f(y)}{\|f(y)\|}.
\label{eq:twenty}
\end{equation}

We observe
\begin{equation}
\left\| \lambda_{1}\frac{x}{\|x\|} + \lambda_{2}\frac{y}{\|y\|}
\right\|^2 = \lambda_{1}^{2} + \lambda_{2}^{2} ~=~
\| (\lambda_{1}, \lambda_{2}) \|^{2}.
\label{eq:twentyone}
\end{equation}
This implies that $(\lambda_{1}, \lambda_{2}) \in D_{2}$ if and only
if (\ref{eq:nineteen}) holds true.
On the basis of this fact, let us define a function
$f^*\hbox{:}\ D_{2} \to \R^2$ by
\begin{equation}
f^* (\lambda) = \mu,  \label{eq:twentytwo}
\end{equation}
where $\lambda = (\lambda_{1}, \lambda_{2}) \in D_{2}$ and
$\mu = (\mu_{1}, \mu_{2})$ obey the relation (\ref{eq:twenty}).
Let $\lambda = (\lambda_{1}, \lambda_{2})$ and
$\lambda' = (\lambda'_{1}, \lambda'_{2})$ belong to $D_{2}$ and
let
\begin{equation}
u = \lambda_{1}\frac{x}{\|x\|} + \lambda_{2}\frac{y}{\|y\|}
\quad\mbox{ and }\quad
u' = \lambda'_{1}\frac{x}{\|x\|} + \lambda'_{2}\frac{y}{\|y\|}.
\label{eq:twentythree}
\end{equation}
Then we have
\begin{equation*}
|\langle u, u' \rangle| =
    |\lambda_{1}\lambda'_{1} + \lambda_{2}\lambda'_{2}| =
    |\langle \lambda, \lambda' \rangle|
\end{equation*}
and
\begin{align*}
|\langle f(u), f(u') \rangle|
     &= \left|\left\langle \mu_{1}\frac{f(x)}{\|f(x)\|} +
                                \mu_{2}\frac{f(y)}{\|f(y)\|},\;
                                \mu'_{1}\frac{f(x)}{\|f(x)\|} +
                                \mu'_{2}\frac{f(y)}{\|f(y)\|}
           \right\rangle\,\right|\\[2pt]
     &= |\mu_{1}\mu'_{1} + \mu_{2}\mu'_{2} |\\[2pt]
     &= |\langle f^* (\lambda), f^* (\lambda') \rangle|.
\end{align*}
Since $f$ satisfies the inequality (\ref{eq:one}) for all
$x,y \in D_{n}$, we obtain by (\ref{eq:twentyone}) and
(\ref{eq:twentythree}) that
\begin{align*}
||\langle f^* (\lambda), f^* (\lambda') \rangle| -
|\langle \lambda, \lambda' \rangle||
     &=  ||\langle f(u), f(u') \rangle| -
              |\langle u, u' \rangle||\\[2pt]
     &\leq  \varphi(u,u') = \phi(\|u\|, \|u'\|)\\[2pt]
     &=    \phi(\|\lambda\|, \|\lambda'\|)
              =: \tilde{\phi}(\lambda, \lambda')
\end{align*}
for all $\lambda, \lambda' \in D_{2}$,
where we understand
$\tilde{\phi}\hbox{:}\ \R^2 \times \R^2 \to [0,\infty)$
as a restriction $\varphi |_{\mbox{\scriptsize $\R^2 \times \R^2$}}$.

According to Corollary 12, we get
\begin{equation}
\| f^* (\lambda) \| ~=~ \| \lambda \|  \label{eq:twentyfour}
\end{equation}
for every $\lambda \in D_{2}$.
If we put $\lambda = (\|x\|, 0)$, then $\lambda \in D_{2}$
in view of the assertion that was verified by (\ref{eq:twentyone}).
For this case, it follows from (\ref{eq:twenty}) and
(\ref{eq:twentytwo}) that
\begin{equation*}
f^* (\lambda) ~=~ (\| f(x) \|, 0).
\end{equation*}
And (\ref{eq:twentyfour}), together with Lemma 6, yields
\[\| x \| = \| \lambda \| = \| f^* (\lambda) \| ~=~ \| f(x) \|\]
for each $x$ in $D_{n}$.\hfill$\Box$
\end{proof}

At last, by making use of Lemmas 9 and 13, and considering Lemma 6,
we can prove the main theorem of this paper.


\begin{theorem}[\!]
If a function $f\!:\!D_{n} \to \R^n$
satisfies the inequality
\begin{equation*}
||\langle f(x), f(y) \rangle| - |\langle x, y \rangle||
\leq \varphi(x,y)
\end{equation*}
for all $x,y \in D_{n}${\rm ,} then $f$ satisfies the generalized
orthogonality equation
\begin{equation*}
|\langle f(x), f(y) \rangle| = |\langle x, y \rangle|,
\end{equation*}
for all $x,y \in D_{n}$.
\end{theorem}

Let $B$ be an open ball in $\R^n$ with radius $d>0$ and
centered at the origin, i.e.,
\begin{equation*}
B := \{ x \in \R^n\hbox{:}\; \| x \| < d \}.
\end{equation*}
In view of Theorem 14, the following corollaries are obvious.


\begin{coro}$\left.\right.$\vspace{.5pc}

\noindent If a function $f\!:\!B \to \R^n$
satisfies the inequality
\begin{equation*}
||\langle f(x), f(y) \rangle| - |\langle x, y \rangle||
    \leq \varepsilon\| x \|^{p}\,\| y \|^{p}
\end{equation*}
for some $\varepsilon \geq 0${\rm ,} $p>1$ and for all $x,y \in B${\rm ,}
then $f$ satisfies the generalized orthogonality equation
$(\ref{eq:zero})$ for all $x,y \in B$.
\end{coro}


\begin{coro}$\left.\right.$\vspace{.5pc}

\noindent If a function
$f\!:\!\R^{n}\!\setminus\! B \to \R^n$
satisfies the inequality
\begin{equation*}
||\langle f(x), f(y) \rangle| - |\langle x, y \rangle||
    \leq \varepsilon\| x \|^{p}\,\| y \|^{p}
\end{equation*}
for some $\varepsilon \geq 0${\rm ,} $p<1$ and for all
$x,y \in \R^{n}\!\setminus\! B${\rm ,}
then $f$ satisfies the generalized orthogonality equation
$(\ref{eq:zero})$ for all $x,y \in \R^{n}\!\setminus\! B$.
\end{coro}

If we assume $p=0$ in Corollary 16, then we can extend the result
of Chmieli\'{n}ski [4] which was introduced in \S1 to the
case of restricted (unbounded) domains.


\section{Applications}

In this section, we will still use the notations $D_n$ and
$\varphi$ to denote the ones defined in \S1.
With these notations, we will prove the superstability
of the orthogonality equation
\begin{equation}
\langle f(x), f(y) \rangle ~=~ \langle x, y \rangle
\label{eq:twentysix}
\end{equation}
on restricted domains.
Every solution of the orthogonality equation (\ref{eq:twentysix})
is an isometry.

We will first improve Lemma 9 adequately for our purpose.


\begin{lemma}
If a function $f\hbox{:}\ D_{n} \to \R^n$
satisfies the inequality
\begin{equation}
|\langle f(x), f(y) \rangle - \langle x, y \rangle| ~\leq~
\varphi(x,y)
\label{eq:twentyseven}
\end{equation}
for all $x,y \in D_{n}${\rm ,} then it holds
$\,\cos A(f(x), f(y)) = \cos A(x, y)$
for any $x, y \in D_{n}\!\setminus\!\{ 0 \}$.
\end{lemma}

\begin{proof}
Since the inequality (\ref{eq:twentyseven}) implies
the validity of the inequality (\ref{eq:one}), all lemmas, theorems
and corollaries in the previous sections hold true for this case.

Let $x \in D_{n}\!\setminus\!\{ 0 \}$ be given.
According to Lemma 4, there exists a function
$\mu_{x}\hbox{:}\ \R \to \R$ with
\begin{equation}
f(c^{-k}x) ~=~ \mu_{x}(c^{-k}) f(x)
\label{eq:twentyeight}
\end{equation}
for any $k \in \N$.
If we replace $x$ and $y$ in (\ref{eq:twentyseven}) by $c^{-k}x$
and $x$, respectively, then it follows from (ii) that
\begin{align*}
|c^{k} \mu_{x}(c^{-k})\| f(x) \|^2 - \|x\|^{2}\,|
     & \leq  c^{k} \varphi(c^{-k}x, x)\\[2pt]
     & \to   0 \;\;\mbox{ as } k \to \infty.
\end{align*}
Since $\| f(x) \|^2 > 0$ and $\|x\|^2 > 0$, we have
\begin{equation}
\mu_{x}(c^{-k}) > 0
\label{eq:twentynine}
\end{equation}
for any sufficiently large $k \in \N$.

By using Lemma 3, Lemma 4, (\ref{eq:twentyeight}) and
(\ref{eq:twentynine}), we get
\begin{equation}
\cos A(f(x), f(y))=\cos A(f(c^{-k}x), f(c^{-k}y))
\label{eq:thirty}
\end{equation}
for all $x, y \in D_{n}\!\setminus\!\{ 0 \}$ and for all sufficiently
large $k \in \N$.

If we replace $x$ and $y$ in (\ref{eq:twentyseven}) by
$c^{-k}x$ and $c^{-k}y$, respectively, and if we multiply
the resulting inequalities by $c^{2k}$, then we obtain
\begin{align*}
&\| x \| \, \| y \| \cos A(x,y)
- c^{2k}\varphi(c^{-k}x, c^{-k}y) \\[2pt]
&\quad \leq c^{2k} \, \| f(c^{-k}x) \| \, \| f(c^{-k}y) \| \,
\cos A(f(c^{-k}x), f(c^{-k}y))\\[2pt]
&\quad \leq \| x \| \, \| y \| \, \cos A(x,y)
+ c^{2k}\varphi(c^{-k}x, c^{-k}y).
\end{align*}

Taking the limit as $k \to \infty$ and using (iii),
(\ref{eq:thirty}) and Lemma 8, we can conclude that
our assertion is valid.\hfill$\Box$
\end{proof}

By using Lemmas 13 and 17 and considering Lemma 6 also,
we will prove the superstability of the orthogonality equation
on restricted domains.


\begin{theorem}[\!]
If a function $f\!:\!D_{n} \to \R^n$
satisfies the inequality
\begin{equation*}
|\langle f(x), f(y) \rangle - \langle x, y \rangle| \leq
\varphi(x,y)
\end{equation*}
for all $x,y \in D_{n}${\rm ,} then $f$ satisfies the orthogonality
equation{\rm ,}
$\langle f(x), f(y) \rangle = \langle x, y \rangle${\rm ,}
for all $x,y \in D_{n}$.
\end{theorem}

Let $B$ be an open ball in $\R^n$ defined by
$B = \{ x \in \R^n\hbox{:}\ \| x \| < d \}$ for a given $d > 0$.
The following corollaries are analogous versions of Corollaries
15 and 16 for the orthogonality equation.


\begin{coro}$\left.\right.$\vspace{.5pc}

\noindent If a function $f\!:\!B \to \R^n$
satisfies the inequality
\begin{equation*}
|\langle f(x), f(y) \rangle - \langle x, y \rangle|
    \leq \varepsilon\,\| x \|^{p}\,\| y \|^{p}
\end{equation*}
for some $\varepsilon \geq 0${\rm ,} $p>1$ and for all $x,y \in B${\rm ,}
then $f$ satisfies the orthogonality equation
$(\ref{eq:twentysix})$ for all $x,y \in B$.
\end{coro}


\begin{coro}$\left.\right.$\vspace{.5pc}

\noindent If a function
$f\!:\!\R^{n}\!\setminus\! B \to \R^n$
satisfies the inequality
\begin{equation*}
|\langle f(x), f(y) \rangle - \langle x, y \rangle|
    \leq \varepsilon\,\| x \|^{p}\,\| y \|^{p}
\end{equation*}
for some $\varepsilon \geq 0${\rm ,} $p<1$ and for all
$x,y \in \R^{n}\!\setminus\! B${\rm ,}
then $f$ satisfies the orthogonality equation
$(\ref{eq:twentysix})$ for all $x,y \in \R^{n}\!\setminus\! B$.
\end{coro}

It will be an interesting problem to investigate what happens
if $p = 1$ in the above Corollary 15, 16, 19 or 20.

\section*{Acnowledgement}

The first author was supported by Korea Research Foundation Grant (KRF-2003-015-C00023).


\begin{thebibliography}{99}

\bibitem{1} Alsina C and Garcia-Roig J L,
               On continuous preservation of norms and areas,
               {\it Aequ. Math.} {\bf 38} (1989) 211--215
\bibitem{2} Alsina C and Garcia-Roig J L,
               On the functional equation
                    $| T(x) \circ T(y) | = | x \circ y |$,
               in: Constantin Carath\'{e}odory: An International Tribute
               (ed.) Th~M~Rassias (Singapore: World Scientific) (1991)
               pp.~47--52

\bibitem{3}  Baker J, Lawrence J and Zorzitto F,
               {The stability of the equation $f(x+y) = f(x)f(y)$},
               {\it Proc. Am. Math. Soc.} {\bf 74} (1979) 242--246
\bibitem{4}  Chmieli\'{n}ski J,
               {On the superstability of the generalized
                    orthogonality equation in Euclidean spaces},
               {\it Ann. Math. Sil.} {\bf 8} (1994) 127--140
\bibitem{5}  Chmieli\'{n}ski J,
               {On the stability of the generalized orthogonality
                    equation},
               in: Stability of mappings of Hyers-Ulam type
               (eds) Th~M~Rassias and J~Tabor
               (Palm Harbor, Florida: Hadronic Press) (1994) pp.~43--57
\bibitem{6}  Hyers D H,
               {On the stability of the linear functional equation},
               {\it Proc. Nat. Acad. Sci. USA} {\bf 27} (1941) 222--224
\bibitem{7}  Lomont J S and Mendelson P,
               {The Wigner unitary-antiunitary theorem},
               {\it Ann. Math.} {\bf 78} (1963) 549--559
\bibitem{8}  Rassias Th M,
               {On the stability of the linear mapping
                    in Banach spaces},
               {\it Proc. Am. Math. Soc.} {\bf 72} (1978) 297--300
\bibitem{9}  R\"{a}tz J,
               {Remarks on Wigner's theorem},
               {\it Aequ. Math.} {\bf 47} (1994) 288--289
\bibitem{10} R\"{a}tz J,
               {On Wigner's theorem: Remarks, complements, comments,
                    and corollaries},
               {\it Aequ. Math.} {\bf 52} (1996) 1--9
\bibitem{11} Skof F,
               {Sull'approssimazione delle applicazioni localmente
                    $\delta$-additive},
               {\it Atti Accad. Sci. Torino Cl. Sci. Fis. Mat. Natur.}
               {\bf 117} (1983) 377--389
\bibitem{12} Ulam S M,
               Problems in Modern Mathematics
               (New York: Wiley) (1964)
\bibitem{13} Wigner E P,
               Gruppentheorie und ihre Anwendungen auf die Quantenmechanik
               der Atomspektren
               (Vieweg: Fried) (1931)
\end{thebibliography}
\end{document}